\documentclass[12pt]{article}

\usepackage{graphicx}
\usepackage{amssymb}
\usepackage{amsmath}
\usepackage[colorlinks,bookmarksopen,bookmarksnumbered,citecolor=red,urlcolor=red]{hyperref}

\usepackage[utf8]{inputenc}

\newtheorem{Remarque}{Remark}
\newtheorem{Def}{Definition}
\newtheorem{prop}{Proposition}
\newtheorem{lem}{Lemma}

\begin{document}



\title{On the dispersion relation for compressible Navier-Stokes Equations}

\author{Saad Benjelloun$^{*}$ and Jean-Michel Ghidaglia$^{*,\S}$\\
$^{*}$ Mohammed VI Polytechnic University (UM6P), \\Modeling Simulation \& Data Analytics, Benguerir, Morocco. \\
$^{\S}$Université Paris-Saclay, ENS Paris-Saclay, CNRS, Centre Borelli, \\91190, Gif-sur-Yvette, France.}

\maketitle


\section{Introduction}\label{S:1:intro}

\indent In this paper we revisit the classical sound dispersion and attenuation  theory due to Stokes \cite{Stokes}, 1845, and Kirchhoff \cite{Kirchhoff}, 1868, for the propagation of sound in non-ideal fluids. In particular we reformulate the analysis due to Fletcher \cite{Fletcher}, 1974, showing conditions for which the sound propagates at the isothermal speed of sound. Also we presents asymptotic developments making precise the physical conditions under which the different dispersion and attenuation formulas apply.\\

The more complex case of two-fluid flow is addressed by Benjelloun and Ghidaglia \cite{BG} to which the reader is referred. In \cite{BG}, it is shown that analytical expressions for the speed of sound depend heavily on the chosen model. These sound speed expressions are compared with experimental values. The consequences for CFD models are also discussed.


\section{Dispersion relations for a single conductive and viscous fluid}\label{Sec4}

In Table \ref{tab4} below we first give the main notations that are used in this paper.\\
\begin{table}[ht]
    \centering
    \begin{tabular}{|c|l||c|l|}\hline
       Symbol  & Quantity  & Symbol  &Quantity  \\\hline \hline
      $\rho$   & Density  &$T$   & Temperature \\\hline
      $\mu$ & Dynamic viscosity  & $\lambda$ & Thermal conductivity \\\hline
      $c$   & Adiabatic speed of sound & $c_T$ & Isothermal speed of sound\\\hline
      $C_v$ & Isochoric heat capacity  &  $C_p$ & Isobaric heat capacity \\\hline
      $\gamma$ & Ratio  $C_p/C_v$  &  $\Gamma$ & Gr\"uneisen parameter \\\hline
      $Pr$ & Prandtl number $\mu\,C_p/\lambda$   & $Kn$  & Knudsen number \eqref{Knudsen} \\\hline
    \end{tabular}
    \caption{Notations for a single fluid model (see \eqref{CT} and \eqref{Cv_etc}).}
    \label{tab4}
\end{table}

Considering the 1D compressible Navier-Stokes system for a divariant fluid:
\begin{equation}
\left\{
\begin{array}{l}
 \frac{\partial \rho}{\partial t}+ u \frac{\partial \rho}{\partial x} + \rho \frac{\partial u}{\partial x}=0,\vspace{0.2cm}\\
 \frac{\partial u }{\partial t}+ u \frac{\partial u}{\partial x}
+\frac{1}{\rho} \frac{\partial p}{\partial x} =\frac{1}{\rho}\frac{\partial}{\partial x}\left(\frac{4\,\mu}{3}\frac{\partial u}{\partial x}\right)\,,\vspace{0.2cm}\\
\frac{\partial s}{\partial t}+ u \frac{\partial s}{\partial x} = \frac{4\,\mu}{3\,\rho\, T} \left(\frac{\partial u}{\partial x}\right)^2 + \frac{1}{\rho\, T}  \frac{\partial}{\partial x}\left(\lambda\,\frac{\partial T}{\partial x}\right),\vspace{0.2cm}\label{systeme_NST}
\end{array}
\right.
\end{equation}
where the dynamic viscosity $\mu$\,, the thermal conductivity $\lambda$ and the pressure $p$ are known function of the density and the temperature. The linearization of \eqref{systeme_NST} around a constant solution $\rho_0, u_0=0, T_0\,$ reads:

\begin{equation}
\left\{
\begin{array}{l}
 \frac{\partial \rho}{\partial t}+ \rho_0 \frac{\partial u}{\partial x}=0,\vspace{0.2cm}\\
 \frac{\partial u }{\partial t}
+\frac{1}{\rho_0} \frac{\partial p}{\partial x} =\frac{4\,\mu_0}{3\,\rho_0}\frac{\partial^2 u}{\partial x^2},\vspace{0.2cm}\\
\frac{\partial s}{\partial t}=\frac{\lambda_0}{\rho_0\,T_0} \frac{\partial^2 T}{\partial x^2} = \frac{\lambda_0\Gamma_0}{\rho_0^2}\frac{\partial^2 \rho}{\partial x^2}+ \frac{\lambda_0}{\rho_0\,C_{v,0}}\frac{\partial^2 s}{\partial x^2} ,\vspace{0.2cm}\\
p= c_0^2 \rho + \rho_0\,\Gamma_0 T_0 \,s\,,
\label{systeme_NST_L}
\end{array}
\right.
\end{equation}
where (for the signification see Table \ref{tab4}):
\begin{equation}\label{Cv_etc}
\Gamma\equiv \frac1\rho \frac{\partial p}{\partial e} \big|_{\rho}\,,\quad c^2 \equiv \frac{\partial p}{\partial \rho} \big|_s\,,\quad
C_{v}\equiv T\,\frac{\partial s}{\partial T} \big|_{\rho}=\frac{\partial e}{\partial T} \big|_{\rho}\,.
\end{equation}
\begin{Remarque}In general the linearization should be made around a constant solution $\rho_0\,,$ $ u_0\,,$ $T_0\,,$ but using the Galilean invariance of Navier-Stokes equations we can assume, as we did, that $u_0=0$\,. Hence in the results found, all the speed of sound in this paper should be shifted by $u_0$.
\end{Remarque}
The second equality in the $s$ evolution equation in \eqref{systeme_NST_L} follows from the two thermodynamic identities \eqref{Relation_8} given in Appendix \ref{AppenD} \cite{SB}.

The linear differential system (\ref{systeme_NST_L}) is of the form: $$\frac{\partial W}{\partial t} + A \frac{\partial W}{\partial x}= B \frac{\partial^2 W}{\partial x ^2}\,,$$ where $W = ^t(\rho\,,u\,, s)$ and $A$ and $B$ are the $3\times 3$ matrices:
$$
A = \left( \begin{array}{ccc}
0 & \rho_0 & 0 \\
\frac{c_0^2}{\rho_0} & 0 & \Gamma_0\,T_0  \\
0 & 0 & 0 \end{array} \right),\quad
B =
\left( \begin{array}{ccc}
0 & 0 & 0 \\
0 & \frac{4\,\mu_0}{3\,\rho_0} & 0  \\
\frac{\lambda_0\Gamma_0}{\rho_0^2} & 0 & \frac{\lambda_0}{\rho_0\,C_{v,0}} \end{array} \right).
$$
We look for non-vanishing plane-wave solutions of the form:
\begin{equation}\label{plane_wave}
    W= W_0\,\exp{i(k\,x-\omega\,t)}=W_0\,\exp{ik\left(x-\frac{\omega}{k}\,t\right)}\,,
\end{equation}
with $k\in\mathbb{R}$, and $\omega = \omega_R + i \omega_I \in\mathbb{C} $ . Hence ${\omega_R}/{k}$ represents the speed of propagation of the plane-wave and the dispersion relation corresponds to the vanishing of the characteristic polynomial of  the $3\times 3$ matrix:
  $$-i\omega I + i k A + k^2 B\,. $$
The dispersion relation reads, dropping the 0 subscript, as:
\begin{eqnarray}\nonumber
    \left(\frac{\omega}{k}\right)^{3} +   \frac{i\,k }{\rho} \left(\frac{4\,\mu}{3} + \frac{\lambda}{C_v}\right)\left(\frac{\omega}{k}\right)^{2}& -&  \left(c^{2} + \frac{4\,k^{2}\, \lambda\, \mu}{3\,\rho^{2}\,C_v }\right)\left(\frac{\omega}{k}\right)+\\
   & -& \frac{i\, k\,  \lambda}{ \rho}\left(\frac{c^2}{C_v}-\Gamma^2T\right)  = 0\,.
   \label{RD_NS}
    \end{eqnarray}
\section{On the speed of sound in single fluid}

Following L. Landau and E. Lifchitz \cite{Landau} a sound wave is an oscillatory motion with small amplitude in a single compressible fluid. Then these Authors derive the classical linear wave equation where the speed of sound is the usual thermodynamic coefficient which holds this name. Extending this definition to the mixture of two non miscible fluids leads us to consider waves in such a medium. According to Whitham \cite{Whitham}, there is no single precise definition of what exactly constitutes a wave. Nevertheless this Author proposes to distinguish between hyperbolic waves, see Definition \ref{def1}, and dispersive waves. The later being plane-waves where frequency $\omega$ is a defined real function of the wave number $k$ and the function $\omega (k)$ is determined by the particular system under consideration. In this case, the speed of the wave is its phase speed, that is $\omega (k)/k$, see (\ref{plane_wave}), and the waves are usually said "dispersive" if this phase speed is not a constant but varies with $k$.

\paragraph{Dispersion relation}
Considering the classical compressible Navier-Stokes equation in $1D$, as written in \eqref{systeme_NST}, we have shown that small disturbances around a constant state of rest will propagates as the superposition of plane waves of the form \ref{plane_wave} ($k\in\mathbb{R}\,,\omega\in\mathbb{C}$)
provided $\omega$ and $k$ satisfy equation 
\eqref{RD_NS}, which is termed as the dispersion relation.
Setting (the notations are given in Table \ref{tab4}):
\begin{equation}\label{CT}
a = \frac{2\,k\, \mu}{3\, \rho}\equiv a_1\,k\quad
\mbox{ and }\quad c_T^2  \equiv\frac{\partial p}{\partial \rho}\big|_T\,,
\end{equation}
where $a$ is Stokes' attenuation as in Stokes \cite{Stokes},
the dispersion relation \eqref{RD_NS} can be written as:
\begin{equation}
\label{NS1D}
\left(\frac{\omega}{k}\right)^{3} + 2\,i\,k\,a_1  \,\left( 1+\frac{3\,\gamma}{4\,Pr} \right)\left(\frac{\omega}{k}\right)^{2}  -  \left(c^{2} + \frac{3\,\gamma}{Pr}a_1^2\,k^2\right)\left(\frac{\omega}{k}\right) - i\,\frac{3\,c^2\,k\,a_1}{2\,Pr} = 0\,,
\end{equation}
thanks also to the identity $c_T^2= c^2 - \Gamma^2 C_v T ={c^2}/{\gamma}$\,, see \eqref{Relation_4} in Section \ref{AppenD}.


\subsection{Euler's equation, the hyperbolic case}
For $\mu= 0$ and $\lambda=0$ (or $a_1=0$, $Pr=\infty$), we recover the speeds of propagation ${\omega}/{k} \in \{ 0, \pm c  \}$ for the inviscid Euler equation, with $c$ the adiabatic speed of sound for the fluid. This system is hyperbolic and non dispersive. Hyperbolicity can be defined as follows.\\

Considering a linear differential system of the form:
\begin{equation}
\label{linearsystem}
\frac{\partial W}{\partial t} + A\, \frac{\partial W}{\partial x} = 0\,,
\end{equation}
where $A$ is a $N\times N$ constant matrix and $N$ is the number of differential equations occurring in the model, we recall the following.
\begin{Def}\label{def1}
The model \eqref{linearsystem} is said to be hyperbolic if there exists a basis $(r_1\,,\ldots\,,r_N)$ of $\mathbb{R}^N$ made with eigenvectors of $A$\,:
\begin{equation}
A\, r_k = \lambda_k\,r_k\,,\quad k=1\,,\ldots, N\,,
\end{equation}
where $\lambda_k\in \mathbb{R}$ are the associated eigenvalues.
\end{Def}
Dimensional analysis shows that the dimension of the $\lambda_k$ is $m/s$
{\it i.e.} they are velocities. Looking for non-vanishing plane-wave solutions, (\ref{plane_wave}), leads to the simple dispersion relation:
\begin{equation}\label{plane_wave1}
    \frac{\omega}{k}\in \{\lambda_1\,,\ldots\,,\lambda_N\}\,,
\end{equation}
{\it i.e.} all these models are non dispersive equations since ${\omega}/{k}$ is constant.


\subsection{Stokes' model: the non-conductive case}
In the case of an non-conductive ($\lambda=0$\,, $Pr=\infty$) viscous flow  ($\mu \neq 0$), we recover Stokes' attenuation and dispersion relations as in Stokes \cite{Stokes}:
 \begin{equation}\label{RDStokes}
 \left(\frac{\omega}{k}\right)^{3}  +  \frac{4\,i\, k\,\mu}{3\,\rho} \left(\frac{\omega}{k}\right)^{2} -  c^{2}\,\left(\frac{\omega}{k}\right) = 0\,,\end{equation}
  \begin{equation}\label{RDStokes1}
  \Rightarrow  \frac{\omega}{k} \in \left\{ 0, {- i \frac{2\,k\, \mu}{3\rho} } \pm  { \sqrt{ c^2 - \frac{4\,k^2 \mu^2}{9\,\rho^2} } }\right\}.    \end{equation}
Let us consider a plane wave \eqref{plane_wave}. Here $1/k$ represents a characteristic length for the wave. In the context of continuum mechanics (by opposition with rarefied flows) in which the Navier-Stokes equations are valid, the Knudsen number, $Kn$, built with this characteristic length should be not greater than the critical Knudsen number $\text{Kn}_c=10^{-2}$\,, that is:
\begin{equation}\label{Knudsen}
    K_n\equiv\frac{k\,\mu}{\rho\,c}\le \text{Kn}_c=10^{-2} \,.
\end{equation}
Hence:
\begin{eqnarray}\nonumber
    \sqrt{ c^2 - \frac{4\,k^2 \mu^2}{9\,\rho^2} }
&=&c\sqrt{ 1 - \frac{4\,(K_n)^2}{9} }\approx c\,\left(1- \frac{2\,(K_n)^2}{9}\right)\\
&=&c\,\left(1- \frac{2\,k^2 \mu^2}{9\,\rho^2\,c^2}\right)=c-\frac{2\,k^2 \mu^2}{9\,\rho^2\,c}\,.
\end{eqnarray}
According to \eqref{plane_wave}, the disturbances are linear combination of functions of the form:
 \begin{equation}\label{plane_wave_NS}
    W= W_0\,\exp{-\frac{2\,k^2\,\mu\,t}{3\,\rho}}\exp{ik\,(x\pm c(k)\,t)}\,,
\end{equation}
\begin{equation}\label{plane_wave_NS0}
    W= W_0\,\exp{-\frac{2\,k^2\,\mu\,t}{3\,\rho}}\exp{ik\,x}\,.
\end{equation}
The waves \eqref{plane_wave_NS} are dispersive since their velocities $\pm c(k)$\,, where:
\begin{equation}\label{c_de_k}c(k)\equiv { \sqrt{ c^2 - \frac{4\,k^2 \mu^2}{9\,\rho^2} } }= \sqrt{ c^2 - k^2\,a_1^2}
\approx c-\frac{2\,k^2 \mu^2}{9\,\rho^2\,c} \,,
\end{equation}
depends on $k$\,. Moreover the waves are exponentially damped with time (and propagation in space), with the characteristic Stokes attenuation $a_1$.

\begin{table}[h!]
\begin{center}
\begin{tabular}{cccccc}
\hline
Property       & Air & Freon &  Water & Honey & Mercury\\
\hline
$\rho$ (Kg/ $m^3$ ) & 1.225 & 1570 & 997 & 1400 & 13500 \\
$\mu$ (Pa.s)  & $1.81 \, 10^{-5}$ & $2,6 10^{-4}$ &  $8.9 \, 10^{-4}$  & 10 & $1.5 \, 10^{-3}$ \\
c (m.$s^{-1}$)  & $340$ & 716  & 1480   & 2030 & 1450\\
$\lambda$ (W/m K) & $2.6 \, 10^{-2}$ & 10 & 0.6 & 0.5  & 8.3 \\
\hline
$\frac{\mu}{\rho c} (m)$ & $4.3\, 10^{-7}$ & $2.3 \, 10^{-10}$ & $6 \, 10^{-10}$ & $4.7 \,10^{-6}$ & $7.6 \, 10^{-11}$\\
\hline
\end{tabular}
\caption{Typical properties at standard conditions for a selection of materials.}
\label{tab2}
\end{center}
\end{table}

\begin{Remarque}
For ideal gases or gases near the ideal state, the Knudsen number defined in \eqref{Knudsen} relates to the molecular free path $\Lambda$, as we have 
$$
K_n = \Lambda\,k \sim \frac{k\mu}{\rho c}.
$$

This is due to the relation $\mu = \rho \,\Lambda\, \sqrt{ \frac{k_b\,T}{2\,\pi\, m}}$, a consequence of the Maxwell-Boltzman distribution, with $m$ the molecular weight of the gaz, and the relation $c =\sqrt{ \frac{\gamma\, k_b\, T}{m}} $. Here we extend the definition \eqref{Knudsen} to all fluids. Expression \eqref{Knudsen} gives and evaluation of $k \times 4.3\, 10^{-7}$ for air and $k \times 6 \, 10^{-10}$ for water at $25 ^o C$. For ultrasound at 20MHz, the order of magnitude for $k$ is $10^{-2}$ to $10^{-3}$  for air and water at standard conditions and $Kn$ is hence very small even for high frequency sounds. See Table \ref{tab2} for an estimation of $\frac{\mu}{\rho\,c}$ in formula \eqref{Knudsen} for other materials.
\end{Remarque}


\subsection{The non viscous case}

For $\mu=0$ and $\lambda\neq 0$ the dispersion relation \eqref{RD_NS} reads:
    \begin{equation}
    \left(\frac{\omega}{k}\right)^{3} +   \frac{i\,k\,\lambda }{\rho\,C_v} \left(\frac{\omega}{k}\right)^{2}-  c^{2} \left(\frac{\omega}{k}\right)-
 \frac{i\, k\,  \lambda\,c^2}{\gamma\, \rho\,C_v} = 0\,.
   \label{GV1}
    \end{equation}
    In Vekstein \cite{Vekstein} this dispersion relation is given in the case of perfect gas. Introducing here a Knudsen number based on thermal diffusivity:
\begin{equation}\label{GV2}
    Kn_{th}\equiv \frac{k\,\lambda}{\rho\,C_p\,c}\,,
\end{equation}
we can rewrite \eqref{GV1} as:
 \begin{equation}\label{GV3}
    \left(\frac{\omega}{k}\right)^{3} +   i\,\gamma\,c\,Kn_{th} \left(\frac{\omega}{k}\right)^{2}-  c^{2} \left(\frac{\omega}{k}\right)-i\,c^3\,Kn_{th} = 0\,.
    \end{equation}
Here again, since we are in the context of fluid dynamics, the Knudsen number is small. For example \eqref{GV2} gives a value of $k \times 6.6\, 10^{-8}$ for air and $k \times 9.6 \, 10^{-11}$ for water at $25 ^o C$. The techniques used in the Sections \ref{AppenC.1} and \ref{AppenC.2} lead immediately to the following asymptotics:
\begin{equation}\label{GV4}
    \frac{\omega^\pm_I}{k}= -\frac{(\gamma-1)\,k\,\lambda}{2\,\rho\,\gamma\,C_v}+\mathcal{O}\left(\frac{k^2\,\lambda^2}{\rho^2\,C_v^2\,c}\right)\,,\quad \frac{\omega^\pm_R}{k}= \pm c+\mathcal{O}\left(\frac{k^2\,\lambda^2}{\rho^2\,C_v^2\,c}\right)\,,
\end{equation}        
        \begin{equation}\label{GV5}
    \frac{\omega_I^0}{k}= -\frac{k\,\lambda}{\rho\,\gamma\,C_v}+\mathcal{O}\left(\frac{k^2\,\lambda^2}{\rho^2\,C_v^2\,c}\right)\,,\quad \frac{\omega_R^0}{k}=\mathcal{O}\left(\frac{k^2\,\lambda^2}{\rho^2\,C_v^2\,c}\right)\,,
\end{equation}
and \eqref{GV4} generalizes Vekstein \cite{Vekstein} to the case of arbitrary divariant fluids. Hence, we see that in this case the speed of propagation, at order zero in $Kn_{th}.c$, is given by the adiabatic speed of sound $c$ and that the attenuation is the same as the one predicted by Stokes-Kirchhoff theory \cite{Kirchhoff}, with vanishing viscosity: 
$$
\frac{(\gamma-1)\,k\,\lambda}{2\,\rho\,\gamma\,C_v} = \frac{(\gamma-1)\,k\,\lambda}{2\,\rho\,C_p} =  \frac{(\gamma-1) c }{2} Kn_{th} = \frac{k\,\lambda}{2\,\rho\,C_p} (\frac{c^2}{c_T^2}-1)  .
$$

Regarding the speed dispersion given by the expressions for $\omega_R / k$, we notice that it is second order w.r.t. $\text{Kn}_{th}$.

\subsection{On the general case}
In this Section our goal is to study the dispersion relation \eqref{RD_NS} that we have rewritten (see \eqref{NS1D}) as $P\left(\frac{\omega}{k}\right)=0$\,, where (we refer to \eqref{CT} and Table \ref{tab4} for the notation):
\begin{equation}\label{A3.1}
P(X)\equiv X^{3} + 2\,i\,k\,a_1  \,\left( 1+\frac{3\,\gamma}{4\,Pr} \right)X^{2}  -  \left(c^{2} + \frac{3\,\gamma}{Pr}a_1^2\,k^2\right)X - i\,\frac{3\,k\,a_1\,c^2}{2\,Pr}\,.
\end{equation}
Introducing the two second-degree polynomials $Q$ and $Q_T$\,:
\begin{eqnarray}\label{A3.2}
Q(X)\equiv X^{2}  + 2\,i\,k\,a_1  \, X -  c^{2}\,,\\
\label{A3.3}
Q_T(X)\equiv X^{2}  + 2\,i\,k\,a_1  \, X -  \frac{c^2}{\gamma}\,,
\end{eqnarray}
we see that (note that these two polynomials are independent of $Pr$):
\begin{equation}\label{A3.4}
P(X) = X\,Q(X)+i\,k\,a_1\,\frac{3\,\gamma}{2\,Pr}\,Q_T(X)\,.
\end{equation}
\subsubsection{Asymptotic for large Prandtl number}\label{AppenC.1}
When the Prandtl number $Pr$ is infinite, that is the case of an non conductive ($\lambda=0$) viscous flow  ($\mu \neq 0$), we have $P(X)=X\,Q(X)$ and we recover \eqref{RDStokes}.\\
Since the $3$ roots of $P$ in this case are distinct, one can see easily\footnote{Using the Implicit Function Theorem as we do hereafter in the proof of Proposition \ref{prop20}.} that for large Prandtl number ({\it i.e.} for $\lambda << C_p\,\mu$), that is the case where viscous effects are much more preponderant w.r.t. thermal ones, the $3$ roots of $P$ depends smoothly on $1/Pr$ and can be expanded as:
\begin{equation}\label{A3.5}
X_{Pr}^{(\ell)}=X_0^{(\ell)}-i\,k\,a_1\,\frac{3\,\gamma}{2\,Pr}
\frac{Q_T(X_0^{(\ell)})}{Q(X_0^{(\ell)})+X_0^{(\ell)}\,Q'(X_0^{(\ell)})}+\mathcal{O}\left(\frac{1}{Pr^2}\right)\,,
\end{equation}
where the $X_0^{(\ell)}$ are the $3$ roots of $X\,Q(X)$. According to \eqref{RDStokes1}, $c(k)$ is defined in \eqref{c_de_k}, these roots are:
\begin{equation}\label{A3.6}
X_0^{(-1)}=- i k\,a_1 -  c(k)\,,\quad X_0^{(0)}=0\,,\quad X_0^{(+1)}=-ik\,a_1 +  c(k)\,,
\end{equation}
and after some computations we deduce from \eqref{A3.5} the following result.
\begin{prop}\label{prop15}
For large Prandtl number, the three roots of $P$ satisfy:
\begin{eqnarray}\label{A3.7}
X_{Pr}^{(\mp)}=-ik\,a_1 \mp  c(k)
\mp\frac{3\,(\gamma-1)\,k\,a_1\,(k\,a_1\pm ic(k))}{4\,c(k)}\frac{1}{Pr}+\mathcal{O}\left(\frac{1}{Pr^2}\right)\,,\\
\label{A3.8}
X_{Pr}^{(0)}=-i\,\frac{3\,k\,a_1}{2\,Pr}+\mathcal{O}\left(\frac{1}{Pr^2}\right)\,.
\end{eqnarray}
\end{prop}
\subsubsection{Asymptotic for small Prandtl number}\label{AppenC.2}
Let us now address the case where the Prandtl number is small ({\it i.e.} for $C_p\,\mu << \lambda$), that is the case where thermal effects are much more preponderant w.r.t. viscous ones.\\
At the limit $Pr=0$, the dispersion relation $P\left(\frac{\omega}{k}\right)=0$ reads $Q_T\left(\frac{\omega}{k}\right)=0$ but $Q_T$ has only two roots while $P$ has three. We are going to prove that for $Pr<<1$ the two roots of $P$ will be on curves starting from the two roots of $Q_T$ (as in \eqref{A3.5}) and the third one is large $\mathcal{O}\left(\frac{1}{Pr}\right)$, see Propositions \ref{prop10} and \ref{prop20}.\\
Indeed denoting by $\xi^{(\ell)}_{Pr}$\,, $\ell=-1\,,0$ and $1$, the three roots of $P$, and if $\xi_{Pr}^{(- 1)}$ (resp. $\xi_{Pr}^{(+1)}$) is close to $\xi^{(- 1)}$ (resp. $\xi^{(+1)}$) where $\xi^{(\ell)}$ are the roots of $Q_T$ that is:
\begin{equation}\label{A3.11}
\xi^{(\ell)}\equiv - i k\,a_1 +\ell \, c_T(k)\,,\quad c_T(k)\equiv \sqrt{ c_T^2 - k^2\,a_1^2}\,,\quad \ell=\pm 1\,,
\end{equation}
we are going to prove the following result.
\begin{prop}\label{prop10}
We have the following asymptotic behavior:
\begin{eqnarray}\label{A9.1}
\lim_{Pr\rightarrow 0}\xi_{Pr}^{(\mp)}=- i k\,a_1 \mp \, c_T(k)\,,\\
\label{A9.2}
\xi_{Pr}^{(0)}\sim i\frac{3\,k\,a_1\,\gamma}{2\,Pr}\,,\mbox{ as } {Pr\rightarrow 0}\,
.
\end{eqnarray}
\end{prop}
This result follows readily from a simple observation.
\begin{lem}
The three roots $\xi^{(\ell)}_{Pr}$ of $P$ satisfy the identity:
\begin{equation}\label{A3.10}
\xi_{Pr}^{(-1)}\,\xi_{Pr}^{(0)}\,\xi_{Pr}^{(+1)}=-i\frac{3\,k\,a_1\,c^2}{2\,Pr}\,.
\end{equation}
\end{lem}
This relation is obvious since according to \eqref{A3.4}, $P(0)=i\,k\,a_1\,\frac{3\,\gamma}{2\,Pr}\,Q_T(0)$ and then Proposition \ref{prop10} follows.\\
Here again, in the spirit of Proposition \ref{prop15}, we can refine Proposition \ref{prop10} and prove the following result.
\begin{prop}\label{prop20}
For small Prandtl number, the three roots of $P$ satisfy:
\begin{eqnarray}\label{A3.12}
\xi_{Pr}^{(\mp)}=- i k\,a_1 \mp \, c_T(k)\pm \frac{(\gamma-1)\,c^2\,(k\,a_1\mp i\,c_T(k))}{3\,\gamma^2\,k\,a_1\,c_T(k)}Pr+\mathcal{O}(Pr^2)\,,\\
\label{A3.13}
\xi_{Pr}^{(0)}= i\frac{3\,k\,a_1\,\gamma}{2\,Pr}+ \mathcal{O}(Pr)\,
.
\end{eqnarray}
\end{prop}
\paragraph{Proof} The $\xi_{Pr}^{(\ell)}$ are the solutions of $P(X)=0$ or equivalently:
\begin{equation}\label{A3.15}
F(X\,,Pr)\equiv Pr\,X\,Q(X)+i\,\frac{3\,k\,a_1\,\gamma}{2}\,Q_T(X)=0\,.
\end{equation}
The function $F$ is smooth and $\frac{\partial F}{\partial X}(X\,,0)=i\,\frac{3\,k\,a_1\,\gamma}{2}\,Q'_T(X)$\,. For each $\ell=\pm 1$, $F(\xi^{(\ell)}\,,0)=i\,\frac{3\,k\,a_1\,\gamma}{2}\,Q_T(\xi^{(\ell)})=0$ and since the roots $\xi^{(\ell)}$ of $Q_T$ are simple: $\frac{\partial F}{\partial X}(\xi^{(\ell)}\,,0)\neq 0$. Hence by the Implicit Function Theorem, for small $Pr$, there exists two smooth curves such that $F(\xi_{Pr}^{(\ell)}\,,Pr)=0$ and $\xi_{0}^{(\ell)}=\xi^{(\ell)}$\,. Then \eqref{A3.12} follows immediately from the first order Taylor expansion of $F$ with respect to $X$ and $Pr$ at the point $(\xi^{(\ell)}\,,0)\,.$

Concerning  \eqref{A3.13}, we simply use \eqref{A3.12} together with the identity \eqref{A3.10}.

 \subsubsection{Asymptotic for the speed of sound}
 The general dispersion relation (\ref{NS1D}) can be explicitly solved using Cardano formula since it is a third order polynomial equation in $\omega/k$ for fixed $k$.
 In general this dispersion relation has $3$ solutions $(\omega\,,k)$ with $k\in\mathbb{R}$ and $\omega\in\mathbb{C}$\,. Writing $\omega=\omega_R+i\,\omega_I$\,, real and imaginary parts, we see that \eqref{plane_wave} reads:
\begin{equation}\label{plane_wave_bis}
    W= W_0\,\exp{(\omega_I\,t)}\,\exp{ik\left(x-\frac{\omega_R}{k}\,t\right)}\,.
\end{equation}
If we stick to the definition by Whitham \cite{Whitham} of dispersive waves, the phase speed of the wave is $\omega_R/k$ while $\omega_I$ corresponds to attenuation (for negative values) or amplification otherwise.
For large or small values of the Prandtl number, Pr, according to Propositions \ref{prop15} and \ref{prop20} in Sections \ref{AppenC.1} and \ref{AppenC.2}, we have the following results for the $3$ roots, $(\omega/k)^\pm$ and $(\omega/k)^0$, of (\ref{NS1D}).
 \begin{itemize}
        \item In the case where viscous effects are much more preponderant w.r.t. thermal ones, {\it i.e.} for $\lambda << C_p\,\mu $\,, or equivalently $Pr>>1$\,,
\begin{equation}\label{Pr_grand_1}
    \frac{\omega^\pm_I}{k}= -k\,a_1\left(1+\frac{3\,(\gamma-1)}{{4\,Pr}}\right)+\mathcal{O}\left(\frac{1}{Pr^2}\right)\,,
\end{equation}
\begin{equation}\label{Pr_grand_2}
\frac{\omega^\pm_R}{k}= \pm c(k)\pm\frac{3\,(\gamma-1)\,k^2\,a_1^2}{4\,c(k)}\frac{1}{Pr}+\mathcal{O}\left(\frac{1}{Pr^2}\right)\,,
\end{equation}
        \begin{equation}\label{Pr_grand_3}
    \frac{\omega_I^0}{k}= -\frac{3\,k\,a_1}{2}\frac{1}{Pr}+\mathcal{O}\left(\frac{1}{Pr^2}\right)\,,
\end{equation}
         \begin{equation}\label{Pr_grand_4}
     \frac{\omega_R^0}{k}=\mathcal{O}\left(\frac{1}{Pr^2}\right)\,.
\end{equation}
       \item In the case where thermal effects are much more preponderant w.r.t. viscous one, {\it i.e.} for $C_p\,\mu << \lambda$\,, or equivalently $Pr<<1$\,,
      \begin{equation}\label{Pr_petit_1}
    \frac{\omega^\pm_I}{k}= - k\,a_1-\frac{(\gamma-1)\,c^2}{3\,\gamma^2\,k\,a_1}{Pr}+\mathcal{O}\left({Pr^2}\right)\,,
\end{equation}
\begin{equation}\label{Pr_petit_2}
\frac{\omega^\pm_R}{k}= \pm c_T(k)\mp\frac{(\gamma-1)\,c^2}{3\,\gamma^2\,c_T(k)}{Pr}+\mathcal{O}\left({Pr^2}\right)\,,
\end{equation}
        \begin{equation}\label{Pr_petit_3}
    \frac{\omega_I^0}{k}=\frac{3\,k\,a_1\gamma}{2\,Pr}+\mathcal{O}\left({Pr}\right)\,,
\end{equation}
         \begin{equation}\label{Pr_petit_4}
     \frac{\omega_R^0}{k}= \mathcal{O}\left({Pr}\right)\,.
\end{equation}
 \end{itemize}

 As far as dispersive waves are concerned, the two expansions of interest are \eqref{Pr_grand_1} and \eqref{Pr_petit_1}. It transpires from these two relations that the speed of sound for large Prandt number is at first order the adiabatic speed of sound of the fluid $\sqrt{\frac{\partial p}{\partial \rho} \big|_s}$ while for small Prandtl number the speed of sound is at first order the isothermal speed of sound of the fluid $\sqrt{\frac{\partial p}{\partial \rho} \big|_T}$. This generalizes the classical analysis in  Fletcher \cite{Fletcher}. Actually, \eqref{Pr_grand_2} and \eqref{Pr_petit_2} give more information,
 recalling that the Knudsen number depends linearly on $k$\,: $Kn={k\,\mu}/{\rho\,c}\,$:
 \begin{eqnarray}\nonumber
    \frac{\omega^+_R}{c\,k}= 1-\frac{2\,(Kn)^2}{9} &+&\frac{\gamma-1}{2\,\gamma} \frac{Kn}{Pr}
    +\frac{2\,\,(\gamma-1)}{9\,\gamma}\frac{Kn^2}{Pr}\\\label{Pr_grand_2.1}  &+&\mathcal{O}\left(\frac{(Kn)^4}{Pr}+\frac{1}{Pr^2}\right)\,,\mbox{ as } Pr\rightarrow \infty\,,
\end{eqnarray}
\begin{eqnarray}\nonumber
     \frac{\omega^+_R}{c_T\,k}= 1-\frac{2\,(Kn)^2}{9} &+&\frac{\gamma-1}{3\,\gamma} {Pr}
    -\frac{2\,(\gamma-1)}{27}{Kn^2}{Pr} \\\label{Pr_petit_2.1}
    &+&\mathcal{O}\left({(Kn)^4}+{Pr^2}\right)\,,\mbox{ as } Pr\rightarrow 0\,.
\end{eqnarray}

\section{Appendix: Some classical thermodynamic identities used in this article}\label{AppenD}

The dispersion relations derived in this paper involves some thermodynamic coefficients that depends on the equations of state. All the coefficients are well referenced in literature but the notations are not universal and furthermore in many references, they are expressed in the case of perfect gas for example in this case the Gr\"uneisen coefficient $\Gamma$ is equal to $\gamma-1$ where $\gamma$ is the ratio of the heat capacity at constant pressure to the heat capacity at constant volume, and this instills sometimes confusion in the obtained results.\\

Considering a divariant substance, we have used the following thermodynamic relations (see \cite{SB}) in this paper:
\begin{equation}\label{Relation_5}
    d\,e=C_v\,d\,T+\frac{\gamma\,\Gamma\,p-(\gamma-1)\,\rho\,c^2}{\gamma\,\Gamma\,\rho^2}d\,\rho\,,
\end{equation}
\begin{equation}\label{Relation_6}
   d\,h=\gamma\,C_v\,d\,T+\frac{\Gamma-\gamma+1}{\Gamma\,\rho}\,d\,p\,,\quad h=e+\frac{p}{\rho}\,,
\end{equation}
\begin{equation}\label{Relation_7}
    d\,p=c^2\,d\,\rho+\rho\,\Gamma\,T\,\,d\,s=\frac{c^2}{\gamma}\,d\,\rho+\rho\,\Gamma\,C_v\,d\,T\,.
\end{equation}
Indeed, starting from
\begin{equation}\label{Relation_1}
    d\,e=C_v\,d\,T+\left(\beta+\frac{p}{\rho^2}\right)d\,\rho\,,
\end{equation}
\begin{equation}\label{Relation_2}
   d\,h=\gamma\,C_v\,d\,T+\left(\alpha+\frac{1}{\rho}\right)d\,p\,,\quad h=e+\frac{p}{\rho}\,,
\end{equation}
\begin{equation}\label{Relation_3}
    d\,p=c^2\,d\,\rho+\rho\,\Gamma\,T\,\,d\,s=c_T^2\,d\,\rho+\epsilon\,d\,T\,,
\end{equation}
where $\alpha$\,, $\beta$\,, \ldots\, can be seen as partial derivatives {\it e.g.} $\epsilon\equiv\frac{\partial p}{\partial T}\big|_\rho$ when the pressure $p$ for this divariant substance is seen as a function of the two independent thermodynamic variables $\rho$ and $T$. We have already identified some of the coefficients in order to be consistent with Table \ref{tab4}, \eqref{CT} and \eqref{Cv_etc}.
Using these two variables, it is elementary to prove the following result.
\begin{prop}\label{prop30}
It follows from Gibbs relation:
\begin{equation}\label{Gibbs}
    T\,d\,s=d\,e-\frac{p}{\rho^2}d\,\rho\,,
\end{equation}
that
\begin{multline}\label{Relation_4}
    \epsilon=\Gamma\,\rho\,C_v \,,\quad \alpha=-\frac{\gamma-1}{\Gamma\,\rho}\,,\quad c_T^2=\frac{c^2}{\gamma}\,,\quad \beta=-\frac{(\gamma-1)\,c^2}{\gamma\,\Gamma\,\rho}\,,\\\Gamma^2\,C_v\,T=\frac{\gamma-1}{\gamma}c^2\,.
\end{multline}
\end{prop}
Hence the identities \eqref{Relation_5} to \eqref{Relation_7}. Combining these identities, we find:
\begin{equation}\label{Relation_8}
d\,p=c^2\,d\,\rho+\rho\,\Gamma\,T\,d\,s\,,\quad d\,T=\frac{\Gamma\,T}{\rho}\,d\,\rho+\frac{T}{C_v}\,d\,s\,,
\end{equation}
which we need to derive \eqref{systeme_NST_L}.




\bibliographystyle{model1-num-names}



\tableofcontents
\end{document}